\documentclass[final,authoryear,12pt]{elsarticle}

\usepackage{amssymb,amsfonts,amsmath}

\usepackage{mathrsfs,bbm}

\usepackage{amsthm}

\newtheorem{definition}{Definition}

\newtheorem{theorem}{Theorem}

\biboptions{round,authoryear}

\usepackage{hyperref}

\hypersetup{colorlinks=true,
    pdftitle={Minimax lower bound for kink location estimators in a nonparametric regression model with long-range dependence},%
    pdfauthor={\textcopyright\ Justin Rory Wishart, The University of Sydney, School of Mathematics and Statistics},%
    pdfsubject={Nonparametric kink estimation},%
    pdfkeywords={nonparametric regression, long-range dependence, kink, minimax},%
    pdfcreator={pdfLaTeX},%
    pdfproducer={LaTeX with hyperref}%
}

\usepackage{fullpage}

\journal{Statistics \& Probability Letters}

\begin{document}

% New commands

\newcommand{\map}[2]{\,{:}\,#1\!\longrightarrow\!#2}
\newcommand{\innerp}[2]{ \langle #1 , #2 \rangle }

\newcommand{\coloneqq}{\mathrel{\mathop:}=}
\newcommand{\eqqcolon}{=\mathrel{\mathop:}}

\begin{frontmatter}



\title{Minimax lower bound for kink location estimators in a nonparametric regression model with long-range dependence}


\author{Justin Rory Wishart}
\ead{justin.wishart@sydney.edu.au}
\ead[http://www.sydney.edu.au/science/maths/u/justinw]{http://www.sydney.edu.au/science/maths/u/justinw}

\address{School of Mathematics \& Statistics F07, The University of Sydney, NSW, 2006, Australia.}

\begin{abstract}
In this paper, a lower bound is determined in the minimax sense for change point estimators of the first derivative of a regression function in the fractional white noise model. Similar minimax results presented previously in the area focus on change points in the derivatives of a regression function in the white noise model or consider estimation of the regression function in the presence of correlated errors. 
\end{abstract}


\end{frontmatter}

% \linenumbers

\section{Introduction}
\label{Intro}

Nonparametric estimation of a kink in a regression function has been considered for Gaussian white noise models by \citet*{Cheng-Raimondo-2008,Goldenshluger-et-al-2008a,Goldenshluger-et-al-2008b}. Recently, this was extended to the fractional Gaussian noise model by \cite{Wishart-2009}. The fractional Gaussian noise model assumes the regression structure,
\begin{equation}
	dY(x) = \mu(x)\,dx + \varepsilon^\alpha dB_H(x), \quad x \in \mathbb{R},
	\label{eq:fixednonparareg}
\end{equation}
where $B_H$ is a fractional Brownian motion (fBm) and $\mu \map{\mathbb{R}}{\mathbb{R}}$ is the regression function. The level of error is controlled by $\varepsilon \asymp n^{-1/2}$ where the relation $a_n \asymp b_n$ means the ratio $a_n/b_n$ is bounded above and below by constants. The level of dependence in the error is controlled by the Hurst parameter $H \in (1/2,1)$ and $\alpha \coloneqq 2 - 2H$, where the i.i.d. model corresponds to $\alpha = 1$. The fractional Gaussian noise model was used by \cite{Johnstone-Silverman-1997,Wishart-2009} among others to model regression problems with long-range dependent errors.

This paper is interested in the performance of estimators of a change-point in the first derivative of $\mu$ observed in model \eqref{eq:fixednonparareg}. This type of change point is called a kink and the location denoted by $\theta$. Let $\widehat \theta_n$ denote an estimator of $\theta$ given $n$ observations. A lower bound is established for the minimax rate of kink location estimation using the quadratic loss in the sense that,
\begin{equation}
 \liminf_{n \to \infty } \inf_{\widehat \theta_n} \sup_{\mu \in \mathscr F_s(\theta)} \rho_n^{-2}\mathbb{E} \left| \widehat \theta_n - \theta\right|^2  \ge C \qquad \text{for some constant $C>0$}.\label{eq:rate}
\end{equation}
The main quantity of interest in this lower bound is the rate, $\rho_n$. In \eqref{eq:rate}, $\inf_{\widehat \theta_n}$ denotes the infimum over all possible estimators of $\theta$. The class of functions under consideration for $\mu$ is denoted $\mathscr F_s(\theta)$ and defined below.
\begin{definition}
	\label{def:functionalclass}
	Let $s \geq 2$ be an integer and $a \in \mathbb{R}\setminus \left\{ 0\right\}$. Then, we say that $\mu\in \mathscr F_s(\theta)$ if, 
	\begin{enumerate}
	\item The function $\mu $ has a kink at $\theta \in (0,1)$. That is, \[\lim_{x \downarrow \theta}\mu^{(1)}(x) - \lim_{x \uparrow \theta}\mu^{(1)}(x) = a \neq 0.\]
	\item The function $\mu \in \mathscr L_2\left(\mathbb{R}\right) \cap \mathscr L_1(\mathbb{R}) $, and satisfies the following condition,
  \begin{equation}
  \int_\mathbb{R} |\widetilde \mu(\omega)||\omega|^s\,d\omega < \infty, \label{sobolev}
\end{equation}
where $\widetilde \mu(\omega) \coloneqq \int_\mathbb{R} e^{-2 \pi i \omega x}\mu(x)\, dx$ is the Fourier transform of $\mu$.
	\end{enumerate}
\end{definition}

The minimax rate for the kink estimators has been discussed in the i.i.d. scenario by \cite{Cheng-Raimondo-2008,Goldenshluger-et-al-2008a} and was shown to be $n^{-s/(2s+1)}$. An extension of the kink estimators to the long-range dependent scenario was considered in \cite{Wishart-2009} that built on the work of \cite{Cheng-Raimondo-2008}. An estimator of kink locations was constructed by \cite{Wishart-2009} and achieved the rate in the probabilistic sense, 
\begin{equation} \left| \widehat \theta_n - \theta\right| = \mathcal O_p (n^{-\alpha s /(2s+\alpha)}),\label{eq:kinkrate}\end{equation}
which includes the result of \cite{Cheng-Raimondo-2008} as a special case with the choice $\alpha = 1$. Both \cite{Cheng-Raimondo-2008} and \cite{Wishart-2009} considered a comparable model in the indirect framework and used the results of \cite*{Goldenshluger-et-al-2006} to infer the minimax optimality of \eqref{eq:kinkrate}. However, the results of \cite{Cheng-Raimondo-2008} and \cite{Wishart-2009} require a slightly more restrictive functional class than $\mathscr F_s(\theta)$. The rate obtained by \cite{Cheng-Raimondo-2008} of $n^{-s/(2s+1)}$ was confirmed as the minimax rate by the work of \cite{Goldenshluger-et-al-2008a} who used the i.i.d. framework and a functional class similar to $\mathscr F_s(\theta)$.

The fBm concept is an extension of Brownian motion that can exhibit dependence among its increments which is typically controlled by the Hurst parameter, $H$ (see \cite*{Beran-1994,Doukhan-et-al-2003} for more detailed treatment on long-range dependence and fBm). The fBm process is defined below.
\begin{definition}
\label{fBm}
	The fractional Brownian motion $\left\{B_H(t) \right\}_{t \in \mathbb{R}}$ is a Gaussian process with mean zero and covariance structure,
	\[
		\mathbb{E} B_H(t)B_H(s) =\frac{1}{2}\left\{ |t|^{2H} + |s|^{2H} - |t-s|^{2H} \right\}. 
	\]
\end{definition}
We assume throughout the paper that $H\in (1/2,1)$, whereby the increments of $B_H$ are positively correlated and are long-range dependent.

 In this paper a lower bound for the minimax convergence rate of kink estimation using the quadratic loss function will be shown explicitly on model \eqref{eq:fixednonparareg}. This is a stronger result in terms of a lower bound than the simple probabilistic result in \eqref{eq:kinkrate} given by \cite{Wishart-2009} and is applicable to a broader class of functions.

\section{Lower bound}
\label{lowerbound}
The aim of the paper is to establish the following result.
\begin{theorem}
\label{thm:lowerboundK}
Suppose $\mu \in \mathscr F_s \left( \theta \right)$ is observed from the model \eqref{eq:fixednonparareg} and  $0 < \alpha < 1$. Then, there exists a positive constant $C < \infty$ that does not depend on $n$ such that the lower rate of convergence for an estimator for the kink location $\theta$ with the square loss is of the form,
 \[ \liminf_{n \to \infty} \inf_{\widehat{\theta}_n} \sup_{\mu \in \mathscr F_s(\theta)}  n^{ 2\alpha s/(2s + \alpha)} \mathbb{E}\left| \widehat{\theta}_n - \theta\right|^2 \ge C. \]
\end{theorem}

From \autoref{thm:lowerboundK} one can see that the minimax rate for kink estimation in the i.i.d. case is recovered with the choice $\alpha = 1$ \citep[see][]{Goldenshluger-et-al-2008a}. Also unsurprisingly, the level of dependence is detrimental to the rate of convergence. For instance as the increments become more correlated, and $\alpha \to 0$, the rate of convergence diminishes.

As will become evident in the proof of \autoref{thm:lowerboundK} the Kullback-Leibler divergence is required between two measures involving modified fractional Brownian motions. To cater for this, some auxiliary definitions to precede the proof of \autoref{thm:lowerboundK} are given in the next section.

\section{Preliminaries}
\label{preliminaries}

In this paper, the functions under consideration are defined in the Fourier domain (see \autoref{def:functionalclass}). Among others, there are two representations for fBm that satisfy \autoref{fBm} that are used in this paper. The first being the moving average representation of \cite{Mandelbrot-van-Ness-1968} in the time domain and second is the spectral representation given by \cite{Samorodnitsky-Taqqu-1994} in the Fourier domain. These both need to be considered since they are both used in the proof of the main result. Both representations have normalisation constants $C_{T,H}$ and $C_{F,H}$ for the time and spectral representations respectively to ensure the fBm satisfies \autoref{fBm}. Start with the time domain representation.
\begin{definition}
\label{fBmMVN}
	The fractional Brownian motion $\left\{B_H(t) \right\}_{t \in \mathbb{R}}$ can be represented by,
	\[
		B_H(t) =\frac{1}{C_{T,H}}\int_\mathbb{R} \left((t-s)_+^{H- 1/2} - (-s)_+^{H-1/2}\right)dB(s),
	\]
 where $C_{T,H} = \Gamma(H + 1/2)/\sqrt{2H \sin (\pi H) \Gamma(2H)}$ and $x_+ = x\mathbbm{1}_{\left\{ x > 0\right\}}(x).$
\end{definition}
%\begin{definition}
%\label{fBmPT}
%The fractional Brownian motion $\left\{B_H(t) \right\}_{t \in [0,1]}$ can be represented by,
%\[
%	B_H(t) =\frac{1}{C_{T,H}}\int_0^t s^{1/2-H}\left(\int_s^t (u-s)^{H - 3/2}u^{H - 1/2}\,du\right)\,dB(s),
%\]
% where $C_{T,H} = \sqrt{\Gamma(2-2H) \sin(\pi(H-1/2))/(2H\pi (H-1/2))}.$
%\end{definition}
For the spectral representation a complex Gaussian measure $\breve B\coloneqq B^{[1]} + i B^{[2]}$ is used where $B^{[1]}$ and $B^{[2]}$ are independent Gaussian measures such that for $i = 1,2;$ $B^{[i]}(A) = B^{[i]}(-A)$ for any Borel set $A$ of finite Lebesgue measure and $\mathbb{E} (B^{[i]}(A))^2 = \text{mesh}(A)/2$. 
\begin{definition}
\label{fBmST}
	The fractional Brownian motion $\left\{B_H(t) \right\}_{t \in \mathbb{R}}$ can be represented by,
	\[
		B_H(t) =\frac{1}{C_{F,H}}\int_\mathbb{R} \frac{e^{i s t} - 1}{is}|s|^{-(H-1/2)}d\breve{B}(s),
	\]
where $C_{F,H} = \sqrt{\pi/(2H \sin (\pi H) \Gamma(2H))}$.
\end{definition}As will become evident in \autoref{proof}, to obtain the lower bound result for the minimax rate, it is crucial to know which functional class to consider for $\mu \map{\mathbb{R}}{\mathbb{R}}$ such that the process $ \int_\mathbb{R} \mu(x)\, dB_H(x)$ is a well defined random variable with finite variance. Two such classes of functions will be considered, $\mathcal H$ and $\widetilde{\mathcal{H}}$, which correspond to the time and spectral versions of fBm respectively. Begin with the moving average representation.
\begin{definition}
 \label{otherstochasticintegralfBmclass}
 Let $H \in \left( 1/2, 1\right)$ be constant. Then the class $\mathcal H$ is defined by,
\[ \mathcal H = \left\{ \mu \map{\mathbb{R}}{\mathbb{R}} \Bigg|  \int_\mathbb{R}\int_\mathbb{R} \mu(x) \mu(y) | x-y|^{-\alpha}\, dy \, dx  < \infty \right\}.\]
 % \innerp{\mu}{\mu}_{\mathcal H} \coloneqq  \frac{(1-\alpha)(2-\alpha)}{2} \int_0^1 \int_0^1 \mu(s)\mu(t) \left| t-s \right|^{-\alpha}\, ds \,dt < \infty \right\}.\]
\end{definition}
Simlar to \eqref{eq:spectralExpectation} there is an inner product on the space $\mathcal H$ that satisfies the following. For all $f,g \in \mathcal H$,
\[ \mathbb{E} \left\{ \int_\mathbb{R} f(x) \, dB_H(x)\int_\mathbb{R} g(y) \, dB_H(y) \right\} = C_\alpha\int_\mathbb{R}\int_\mathbb{R} f(x) g(y)| x-y|^{-\alpha}\, dy \, dx \eqqcolon \innerp{f}{g}_{\mathcal{H}},\]
where the constant $C_\alpha  = \tfrac{1}{2} (1-\alpha)(2-\alpha)$. The other functional class for the spectral representation is denoted by $\mathcal H$ and defined below.
 \begin{definition}
 \label{stochasticintegralfBmclass}
 Let $H \in \left( 1/2, 1\right)$ be constant. Then the class $\widetilde{\mathcal H}$ is defined by,
\[ \widetilde{\mathcal H} = \left\{ \mu \map{\mathbb{R}}{\mathbb{R}} \Bigg| \int_\mathbb{R} |\widetilde\mu(\omega)|^2 \left| \omega \right|^{-(1-\alpha)}\, d\omega < \infty \right\}.\] % \innerp{\mu}{\mu}_{\mathcal H} \coloneqq  \frac{(1-\alpha)(2-\alpha)}{2} \int_0^1 \int_0^1 \mu(s)\mu(t) \left| t-s \right|^{-\alpha}\, ds \,dt < \infty \right\}.\]
\end{definition}
On the space $\widetilde{\mathcal H}$, the stochastic integrals with respect to fBm are well defined and satisfy the following. For all $f,g \in \widetilde{\mathcal{H}}$, 
\begin{equation}
	\mathbb{E} \left\{ \int_\mathbb{R} f(x) \, dB_H(x)\int_\mathbb{R} g(y) \, dB_H(y) \right\} = \frac{1}{C_{F,H}^2}\int_\mathbb{R} \widetilde{f}(\omega) \overline{\widetilde{g}(\omega) }| \omega|^{-(1-\alpha)}\, d\omega \eqqcolon \innerp{f}{g}_{\widetilde{\mathcal{H}}},
	\label{eq:spectralExpectation}
\end{equation}
where $\overline{\widetilde{g}}$ denotes the complex conjugate of $\widetilde{g}$.

These two classes of integrands were considered extensively in \cite{Pipiras-Taqqu-2000}. In this context of this paper the inner products can be used interchangeably because if $\mu \in \mathscr F_s(\theta)$ then $\mu \in \mathscr L_1(\mathbb{R})\cap \mathscr L_2(\mathbb{R})$ and by \citet[Proposition 3.1]{Pipiras-Taqqu-2000} then $\mu \in \mathcal H$. Also, using  \citet[Proposition 3.2]{Pipiras-Taqqu-2000} with the isometry \citet[Lemma 3.1.2]{Biagini-et-al-2008} and Parseval's Theorem then $\mu \in \widetilde{\mathcal H}$ and consequently $\mu \in \mathcal H \cap \widetilde{\mathcal H}$.
\section{Proof of Theorem 1}
\label{proof}
The lower bound for the minimax rate is constructed by adapting the results of \cite{Goldenshluger-et-al-2006} to our framework. This requires obtaining the Kullback-Leibler divergence of two suitably chosen functions $\mu_0$ and $\mu_1$ from the functional class $\mathscr F_s(\theta)$.
 The main hurdle in determining the Kullback-Leibler divergence is the long-range dependent structure in the fBm increments. A summary of  Girsanov type theorems for fBm have been established by \citet*[Theorem 3.2.4]{Biagini-et-al-2008}. Here however, the Radon-Nikodym derivative is the main focus. Once that is determined, the Kullback-Leibler divergence is linked to the lower rate of convergence using \citet[Theorem 2.2 (iii)]{Tsybakov-2009}. Lastly, before proceeding to the proof, the quantity $C > 0$ denotes a generic constant that could possibly change from line to line.

Without loss of generality, consider a function $\mu_0 \in \mathscr F_s(\theta_0)$ where $\theta_0 \in (0,1/2]$ and define $\theta_1 = \theta_0 + \delta$ where $\delta \in (0, 1/2)$ (a symmetric argument can be setup to accommodate the case when $\theta_0 \in [1/2,1)$). Define the functions $v \map{\mathbb{R}}{\mathbb{R}}$ and $v_N \map{\mathbb{R}}{\mathbb{R}}$ such that
\begin{align*}
 v(x) &\coloneqq a((\theta_1\wedge x) - \theta_0) \mathbbm{1}_{ (\theta_0,1] }(x), \qquad v_N(x) \coloneqq \int_{-N}^N \widetilde{v}(\omega) e^{2 \pi i x \omega}\, d\omega,
\end{align*}
where $a$ is the size of the jump given in \autoref{def:functionalclass} and $\widetilde v$ is the Fourier transform of $v$. Note that, $v_N(x)$ is close to $v(x)$ in the sense that it is the inverse Fourier transform of $\widetilde{v}(\omega)\mathbbm{1}_{|\omega| \le N}$ and $\widetilde v_N(\omega) = \widetilde v(\omega)\mathbbm{1}_{|\omega| \le N}$. With these definitions, the derivative takes the form, $ v^{(1)}(x) = a \mathbbm{1}_{[\theta_0, \theta_1] }(x)$ and the function $(\mu_0 - v)$ has a single kink at $\theta_1$. Then define $\mu_1 \coloneqq \mu_0 - (v - v_N)$. The function $v_N$ is infinitely differentiable across the whole real line and smooth for finite $N$, which implies that $\mu_1 = \mu_0 - (v - v_N)$ has a single kink at $\theta_1$. It can be shown that, 
\begin{equation}
	\left|\widetilde{v}(\omega)\right| \le a\delta/(2 \pi  \left|\omega\right|)^{-1}.\label{eq:vomegabound}
\end{equation}
Further, if $N$ is chosen to be $N = \left( s\pi C/(a \delta) \right)^{1/s}$ then $\int_\mathbb{R} | \widetilde{v_N}(\omega)| | \omega|^s \,d\omega < \infty$ and consequently $\mu_1 \in \mathscr F_s(\theta_1)$.
 
To be able to determine the Radon-Nikodym derivative, define $\Delta \coloneqq \mu_0 - \mu_1 = v - v_N$ and note that $\Delta \map{\mathbb{R}}{\mathbb{R}}$. The Radon-Nikodym derivative also needs a paired function
$\underline{\Delta}\map{\mathbb{R}}{\mathbb{R}}$. Define such a function with a singular integral operator with
\begin{equation} \Delta(x) \coloneqq  \varepsilon^{-\alpha}C_\alpha \int_\mathbb{R} |x-y|^{-\alpha} \underline{\Delta}(y) \, dy = \frac{\Gamma(3-\alpha)}{2} \left( \mathcal D_-^{-(1-\alpha)}\underline\Delta(x) + \mathcal D_+^{-(1-\alpha)}\underline\Delta(x)\right),\label{eq:DeltaDecomp}\end{equation}
where, for $\nu \in (0,1) $, $\mathcal D_-^{-\nu}$ and $\mathcal D_+^{-\nu}$ are the left and right fractional Liouville integral operators defined by,
\[ \mathcal D_-^{-\nu} f(x)  \coloneqq \frac{1}{\Gamma(\nu)}\int_{-\infty}^x (x-y)^{\nu -1}f(y)\, dy \qquad \mathcal D_+^{-\nu} f(x)  \coloneqq \frac{1}{\Gamma(\nu)}\int_x^\infty (y-x)^{\nu -1}f(y)\, dy \]
%\[ \widetilde{\underline{\Delta}}(\omega) \coloneqq C_{F,H}^2\varepsilon^{-\alpha} |\omega|^{1-\alpha}\widetilde{\Delta}(\omega). \]
This function $\underline \Delta$ has representation in the Fourier domain with
\[ \widetilde{\underline \Delta} (\omega)  \asymp \varepsilon^{-\alpha}|\omega|^{1-\alpha} \widetilde{\Delta}(\omega).\]
Furthermore, $\underline{\Delta} \in \mathcal H \cap \widetilde{\mathcal H}$. Indeed, by definition, $\Delta = \mu_0 - \mu_1$ with $\mu_0 \in \mathscr F_s(\theta_0)$ and $\mu_1 \in \mathscr F_s(\theta_1)$ which implies that $\Delta \in \mathscr L_1(\mathbb{R}) \cap \mathscr L_2(\mathbb{R})$ and $\widetilde \Delta (\omega) = o(\omega^{-s})$ due to \eqref{sobolev}. First, it will be shown that, $\underline \Delta \in \widetilde{\mathcal H}$.
\begin{align}
	\innerp{\underline \Delta}{\underline \Delta}_{\widetilde{\mathcal H}} &\asymp  \int_\mathbb{R} | \widetilde \Delta (\omega)|^2 |\omega|^{1-\alpha}\, d\omega\nonumber\\
	&\le  C\left\{ \|\Delta\|_{1}^2\int_{|\omega|\le 1} |\omega|^{1-\alpha}\, d\omega +   \int_{|\omega|\ge 1} |\widetilde \Delta (\omega)|^2 |\omega|^{1-\alpha}\, d\omega\right\},\label{eq:uDeltainnerp}
\end{align} 
where $C > 0$ is some constant and $ \|\Delta\|_{1} = \int_\mathbb{R} | \Delta (x)| \, dx $. In \eqref{eq:uDeltainnerp}, the first integral is finite since $\alpha \in (0,1)$ and the last integral is finite since $\widetilde \Delta (\omega) = o(\omega^{-s}) $ for $s \ge 2$, proving $\underline \Delta \in \widetilde{\mathcal H}$. Then apply the isometry in \citet[Lemma 3.1.2]{Biagini-et-al-2008} with Plancherel and \eqref{eq:uDeltainnerp}, it follows that $\underline \Delta \in \mathcal H$.

Now let $P_0$ and $P_1$ be the probability measures associated with model \eqref{eq:fixednonparareg} with $\mu = \mu_0$ and $ \mu = \mu_1$ respectively. 
%That is, $P_0$ is the measure associated with, 
%\begin{align*}
% dY_0(x) &= \mu_0 (x)\, dx + \varepsilon^\alpha \, dB_H(x)
%\end{align*}
%and $P_1$ is the measure associated with, 
%\begin{align*}
% dY_1(x) &= \mu_1 (x) \, dx + \varepsilon^\alpha \, dB_H(x).
%\end{align*}
Define, $\mathring{B}_H(x) \coloneqq \varepsilon^{-\alpha} \int_0^x \Delta(x)\, dx + B_H(x)$. Then under the $P_0$ measure,
\begin{align*}
  dY_0(x) &= \mu_0(x)\,dx + \varepsilon^\alpha \, dB_H(x)
%\\
%  Y_0(x) &= \mathcal{D}^{-1}\mu_0(x) + \varepsilon^\alpha B_H(x)\\
%    &= \mathcal{D}^{-1}\mu_1(x) + \varepsilon^\alpha \left( \frac{\mathcal{D}^{-1}\Delta(x)}{\varepsilon^\alpha} + B_H(x) \right)\\
%    &= \mathcal{D}^{-1}\mu_1(x) + \varepsilon^\alpha \mathring{B}_H(x)\\
= \mu_1(x)\,dx + \varepsilon^\alpha \, d\mathring{B}_H(x).
\end{align*}
%Thus the model $Y_0$ can be modified to be the same as the model $Y_1$ if the measure is used under the $\mathring{B}_H(x)$ process. %Apply \citet[Theorem 3.2.4]{Biagini-et-al-2008} and it follows that the Radon-Nikodym derivative is,
The Radon-Nikodym derivative between these measures takes the form,
\begin{align}
 \frac{d P_1}{d P_0} &\coloneqq \exp \left\{ - \int_\mathbb{R} \underline{\Delta}(x) \, dB_H(x) - \frac{1}{2}\mathbb{E}_{P_0} \left(  \int_\mathbb{R} \underline{\Delta}(x) \, dB_H(x) \right)^2\right\} . \label{eq:radnikderiv}
\end{align}
Indeed to show \eqref{eq:radnikderiv} is valid, for $\underline{\Delta} \in \mathcal H$ and $\psi \in \mathcal H$, use \eqref{eq:DeltaDecomp} and apply \citet[Lemma 3.2.1]{Biagini-et-al-2008} with the change of measure formula in \eqref{eq:radnikderiv} to yield,
\begin{equation}
	\mathbb{E}_{P_1} \left[ \psi(\mathring{B}_H(x)) \right]= \mathbb{E}_{P_0} \left[\psi(\mathring{B}_H(x))\frac{d P_1}{d P_0} \right] = \mathbb{E}_{P_0} \Big[\psi(B_H(x))\Big].
\end{equation}
So, using \eqref{eq:spectralExpectation} in \eqref{eq:radnikderiv}, the Kullback-Leibler divergence between the two models can be evaluated,
\begin{equation}
\mathcal K(P_0,P_1) \coloneqq \mathbb{E} \ln \frac{d P_0}{d P_1} = \frac{1}{2 }  \innerp{\underline{\Delta}}{\underline{\Delta}}_{\widetilde{\mathcal H}}.\label{logstochasticexp}
\end{equation}
To evaluate \eqref{logstochasticexp}, obtain a finer bound on $|\widetilde{\underline{\Delta}}(\omega)|^2$  by recalling that $\Delta = v - v_N$ and using \eqref{eq:vomegabound},
\begin{equation}
  |\widetilde{\underline{\Delta}}(\omega)|^2 \asymp \varepsilon^{-2\alpha} |\widetilde{v}(\omega)|^2 \mathbbm{1}_{\left\{ |\omega| \geq N \right\}} |\omega|^{2 - 2\alpha} \leq \frac{C^2a^2\delta^2}{4\pi^2} \varepsilon^{-2\alpha} |\omega|^{-2\alpha}    \mathbbm{1}_{\left\{ |\omega| \geq N \right\}}. \label{eq:DeltaModulus}
\end{equation}
Apply the bound in \eqref{eq:DeltaModulus} to \eqref{logstochasticexp} with the chosen $N = \left( s\pi C/(a \delta) \right)^{1/s}$,
\begin{align*}
 \mathcal{K}(P_0,P_1)   &= \frac{1}{2} \int_\mathbb{R} |\widetilde{\underline{\Delta}}(\omega)|^2|\omega|^{-(1-\alpha)}\, d\omega\\
  &\leq \frac{ Ca^2 \delta^2}{4 \pi^2} \varepsilon^{-2\alpha}  \int_{\left|\omega\right| \geq N } \left|\omega\right|^{-\alpha-1}\, d\omega\\
%  &= \frac{a^2 \delta^2}{2 \pi^2}  \varepsilon^{-2\alpha} N^{-\alpha}\\
  &= Ca^2 \delta^2\varepsilon^{-2\alpha} \left( s /(a \delta) \right)^{-\alpha/s}\\
  &\asymp \delta^{(2s+\alpha)/s} \varepsilon^{-2\alpha} .
\end{align*}
Now choose $\delta \asymp \varepsilon^{2\alpha s/(2s + \alpha)}$ which guarantees that $\mathcal{K}(P_0,P_1) \leq K < \infty$ for some finite positive constant $K$. Then by \citet[Theorem 2.2 (iii)]{Tsybakov-2009} combined with the fact that $\varepsilon \asymp n^{-1/2}$ it follows that the lower rate of convergence for the minimax risk is $\varepsilon^{2\alpha s/(2s+\alpha)}\asymp n^{- \alpha s/(2s+\alpha)}$. \hfill $\Box$
\section*{Acknowledgements}
The author would like to thank the editor and an anonymous referee for their comments and suggestions which lead to an improved version of this paper.

\end{document}